\documentclass[onecolumn]{paul}

\usepackage{latexsym,amssymb,lastpage}
\usepackage{graphicx,amsfonts}
\usepackage{times,mathptmx,bm,amsmath}
\usepackage{dcolumn}

\newcommand{\rhoeq}{ \rho_{\mbox{\scriptsize eq}} }
\newcommand{\rhoeqone}{ \rho_{\mbox{\scriptsize eq},1} }
\newcommand{\rhoeqtwo}{ \rho_{\mbox{\scriptsize eq},2} }
\newcommand{\rhoeqi}{ \rho_{\mbox{\scriptsize eq},i} }
\newcommand{\phimin}{\phi_{\mbox{\scriptsize min}}}
\newcommand{\lambdamax}{\lambda_{\mbox{\scriptsize max}}}

\renewcommand{\theequation}{\thesection.\arabic{equation}}
\numberwithin{equation}{section}

\begin{document}

\begin{center}
\title{A Paradox of State-Dependent Diffusion\\ and How to Resolve It}
\end{center}

\begin{center}
{{\sc P. F. Tupper and Xin Yang}\\
Department of Mathematics, Simon Fraser University,\\
Burnaby BC, V5A 1S6 Canada.}
\end{center}

\pagestyle{headings}

\maketitle

\begin{abstract}Consider a  particle diffusing in a confined volume which is divided
into two equal regions.  In one region the diffusion coefficient is
twice the value of the diffusion coefficient in the other region.  Will the particle spend equal proportions 
of time in the two regions in the long term?  Statistical
mechanics would suggest yes, since the number of accessible states in
each region is presumably the same.  However, another line of reasoning suggests
that the particle should spend less time in the region with faster
diffusion, since it will exit that region more quickly.  We demonstrate with
a simple microscopic model system that both predictions are consistent
with the information given.  Thus, specifying the diffusion rate as a function of position is not enough to characterize the behaviour of a system, even assuming the absence of external forces.   
We propose an alternative framework for modelling diffusive dynamics in which both the diffusion rate and equilibrium probability density for the position of the  particle are specified by the modeller.
We introduce a numerical method for simulating dynamics in our framework that samples from the equilibrium probability density exactly and is suitable for discontinuous diffusion coefficients. 
\end{abstract}

\section{Introduction} \label{sec:intro}

\begin{figure}
\begin{center}
\includegraphics[width=13cm]{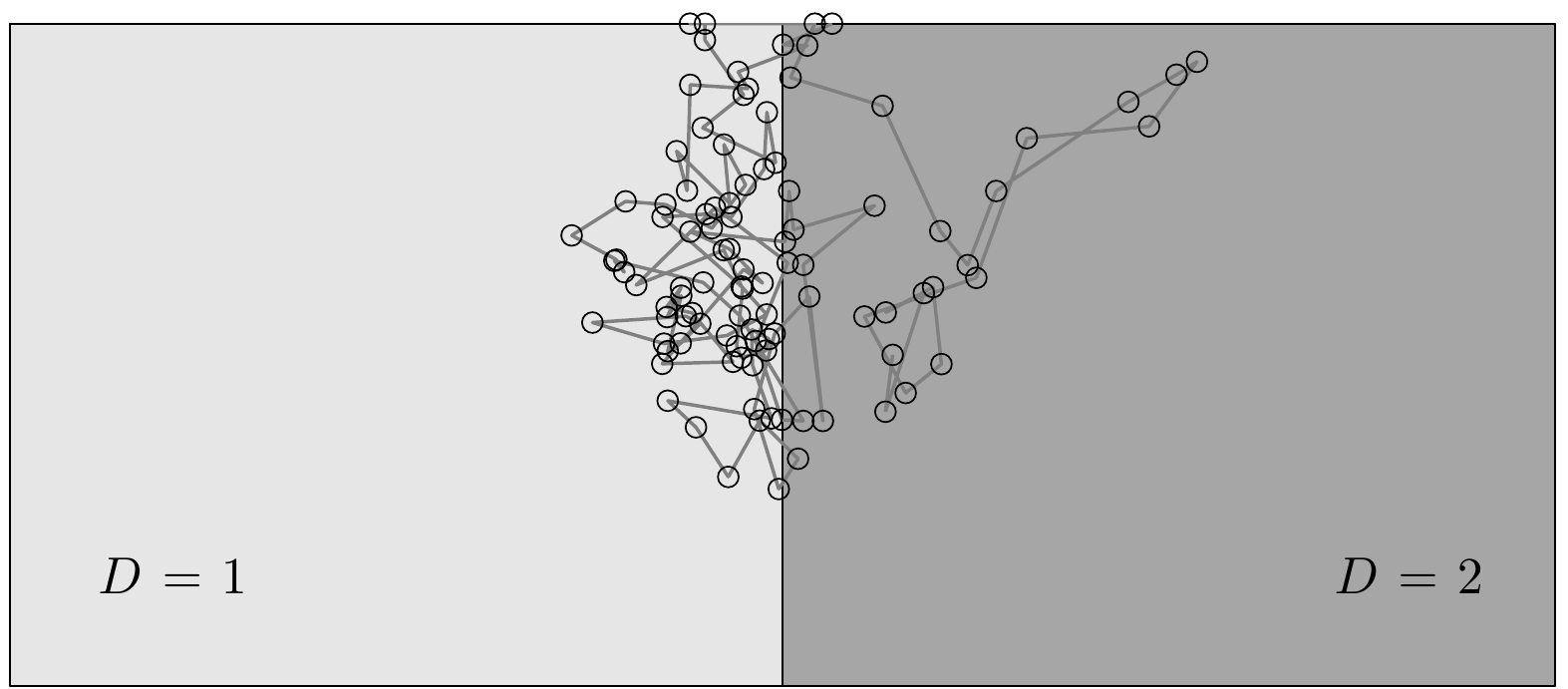}
\end{center}
\caption{  \label{fig:rect_box} Simulation of particle diffusing within a box with reflecting boundary conditions. Points are generated by the Euler-Maruyama method $X_{n+1}= X_n + \sqrt{ 2 h D(X_n) } N_n$, where $h=0.5$, $N_n$ are independent standard two-dimensional Gaussians, and  $D(x)$ is the state-dependent diffusion coefficient. Here  $D(x)=1$ on the left side of the box and $D(x)=2$ on the right side of the box. }
\end{figure}

Consider a particle diffusing in a two-dimensional box with reflecting boundary conditions.
We show a portion of a simulated trajectory of such a system in Figure~\ref{fig:rect_box}.
Suppose that in the left half of the box the particle diffuses with coefficient $D_1$ and that in the right half of the box the particle diffuses with coefficient $D_2 = 2 D_1$.  
We assume that there are no external forces acting on the particle.
Our question is: Does the particle spend an equal fraction of time on each side of the box in the long run?

One answer is based on statistical mechanics.

\begin{center}
\parbox{13cm}{
{\bf Statistical Mechanics Prediction:} The particle will spend an equal proportion of time on each side of the box. 
}
\end{center}

 The justification for this prediction is the principle of statistical mechanics which states that ``An isolated system in equilibrium is equally likely to be in any of its accessible states" (Reif 1965, p. 54).
Since all states in the box are accessible, and there are an equal number of states on each side of the box, the particle should spend an equal proportion of its time on each side of the box.

Another answer is based on the idea of rescaling time in one side of the box.

\begin{center}
\parbox{13cm}{
{\bf Time-Change Prediction:}
The particle will spend less time on the side of the box where the diffusion coefficient is greater.  
}
\end{center}

The justification for this prediction is that, in the absence of any drift, faster diffusion is equivalent to time passing more quickly.
  This means that the periods of time the particle spends on the right side of the box will be shorter than those spent on the left side.  Hence the total time the particle spends on the right hand side of the box will be less.
We show how this prediction is a straightforward consequence of interpreting the particle's motion as drift-free diffusion, where we interpret the state-dependent diffusion coefficient using the It\^{o} convention (Gardiner, 2004).
We can write the equation for the particle's motion as
\begin{equation} \label{eqn:orig_sde}
dx 
= b(x) 
 dB(t),
\end{equation}
where $x=(x_1,x_2)$ and $B$ is standard two-dimensional Brownian motion.
(Equivalently we may write this equation as 
\(
dx/dt = b(x) \eta(t)
\)
where $\eta$ is two-dimensional Gaussian white noise.) 
We specify  $b(x) = b_1$ for $x_1 < 0$ and $b(x) = b_2$ for $x_1>0$. Here $b_i = \sqrt{2D_i}$, $i=1,2$ where $D_i$ is the corresponding diffusion coefficient.
 We enforce reflecting boundary conditions at the four walls of the box.
 The (It\^{o}-)Fokker-Planck equation for $\rho(x,t)$, the probability density of the particle's position at time $t$, is (Gardiner 2004, p. 118)
\[
\frac{\partial}{\partial t} \rho(x,t)
=
\frac{1}{2} \nabla \cdot [ \nabla  \left( b^2(x) \rho(x,t) \right) ]
=  \nabla \cdot [ \nabla \left(  D(x) \rho(x,t) \right) ].
\]
The equilibrium density $\rhoeq(x)$ satisfies
\[
\nabla \left( D(x) \rhoeq(x)\right) = \mbox{const}.
\]
Reflecting boundary conditions for the diffusion correspond to Neumann boundary conditions (zero-flux) for the Fokker-Planck equation.
With these boundary conditions, the unique equilibrium density is
\(
\rhoeq(x) = C/D(x)
\)
for some constant $C$.  Thus, since $D_2=2D_1$, the particle spends half as much time on the right side of the box as on the left.  (We discuss the relation of our question to other interpretations of \eqref{eqn:orig_sde} in Section~\ref{sec:conclusion}.)

Neither the Statistical Mechanics Prediction nor the Time-Change Prediction are definitive.
The Statistical Mechanics Prediction relies on the principle of equal probability of all accessible states, which 
needs to be independently justified for the mesoscopic level of description we are considering here.
 The 
 Time-Change Prediction
  is not definitive since 
   It\^{o} stochastic differential equations  are themselves mesoscopic models  whose use can only be rigorously justified by showing how they arise as the coarse-scale limit of  microscopic dynamics.  Since the two predictions contradict each other, at least one of them must be wrong for any given physical system. The approach by which we resolve this apparent contradiction is to study a microscopic  model  of the box system we have described above and then see what proportion of the time the particle spends on each side in simulations of that system.  If one prediction turned out to always be true for our model system, we could use our result as a baseline for investigating under which more general situations the prediction was still true. However, we will see that even for our simple system there are parameter choices that make the Statistical Mechanics Prediction correct and parameter choices that make the 
   Time-Change Prediction correct.
   This demonstrates that there is no \emph{a priori} reason to conclude that one prediction or another is correct, given only the diffusion coefficient $D(x)$ for the system.

In Section~\ref{sec:model_sys}, we describe our  model system, perform numerical experiments on it, and show that the result of the numerical experiments can be determined analytically from the properties of the model system.
 The system we consider is a deterministic Hamiltonian billiard system: the random Lorentz gas (Dettmann 2000).  The system has two free parameters: disc radius and free volume fraction. 
 There is a one-parameter family of values of these parameters  that can generate the same effective diffusion coefficient.
 We will show how to choose these parameters to obtain arbitrarily good approximations to diffusive motion for the  particle on each side of the box.    Using the degree of freedom in the choice of parameters, we show that the equilibrium density of the  particle is underdetermined by the diffusion coefficient on each side of the region.  For a given $D_1$ and $D_2$, exploiting the flexibility in the parameter choice allows us to create systems in which either the Statistical Mechanics Prediction or the 
 Time-Change
 Prediction is correct. 

Our justification for studying a particular microscopic system is threefold.  
Firstly, there is a solid analytical understanding of the random Lorentz gas on which we can base our simulations.  Secondly, our purpose is not to conclude that a particular style of mesoscopic modelling is always the correct one, but to show that postulating a state-dependent diffusion rate is not enough to fully specify mesoscopic behaviour.  For this objective, it is enough to show that multiple mesoscopic behaviours are possible for a single class of simple models, as we do in this paper.  
Finally, we have resorted to a purely deterministic microscopic model, rather than a  stochastic microscopic model (such as Langevin dynamics or a random walk) to avoid concerns that the method used to introduce randomness at the microscopic level somehow biases the results at the mesoscopic level. 
This last point distinguishes our work from similar discussions of van Kampen (2007, Sec.\ XI.3) and Othmer and Stevens (1997, Sec.\ 2). These authors show that spatially inhomogeneous random walks, with natural choices of parameters, can lead to either prediction holding true at the mesoscopic level. 

Similarly,
the conclusions of our Section~\ref{sec:model_sys} closely parallel those of Korabel and Barkai (2011), in which the same question is answered using a random-walk model on a one-dimensional lattice.  There, following the earlier work of Ovaskainen and Cornell (2003), they show that choosing  how the random walk behaves at the interface of two regions of differing diffusion coefficient leads to different equilibrium probability densities for the system.  Korabel and Barkai explain how to determine the correct interface behaviour of the random walk model using experimentally measurable quantities.  Our approach differs in that, in our model, the interface behaviour is determined indirectly through the microscopic dynamics that we describe.

The results in Section~\ref{sec:model_sys} show that choosing a diffusion coefficient $D(x)$ is not enough to fully specify  diffusive dynamics in the absence of other assumptions. 
In particular, our results show that fixing $D(x)$ is not enough to specify  $\rhoeq(x)$, the equilibrium probability density for the position of the particle.
In Section~\ref{sec:prop} we propose a framework for modelling state-dependent diffusion that makes this fact explicit. Rather than simply specifying a state-dependent diffusion coefficient, we specify diffusion coefficient $D(x)$ and an equilibrium density $\rhoeq(x)$ which together with a detailed balance assumption (no-flux in equilibrium)  completely determine  the dynamics. We then introduce a new method for the numerical simulation of diffusive dynamics which makes use of our framework: the numerical method is expressed in terms of $D(x)$ and $\rhoeq(x)$ and makes no reference to a drift term.
 The method consists of  Euler-Maruyama steps for a purely diffusive It\^{o} stochastic differential equation  together with Metropolis rejections. The method is similar to the Metropolis-adjusted Langevin algorithm 
(Roberts and Tweedie,  1996; Bou-Rabee and Vanden Eijnden, 2010),
 except that there is no drift term in the Euler-Maruyama step and it is the Metropolis rejections that induce any drift in the trajectories.
 The advantages of our method over the Euler-Maruyama method are that it samples space with the correct equilibrium density $\rhoeq(x)$ and performs well even when $D(x)$ and $\rhoeq(x)$ have discontinuities.

One area in which the proposed framework in Section~\ref{sec:prop} could be used is cellular biology. Although earlier models of chemistry in the cytoplasm of cells assumed that chemical species were `well-mixed' and thus ignored diffusion, more recent models have taken the geometry of the cell and the diffusion coefficient of various molecules into account (Turner, Schnell, and Burrage, 2004). Effective diffusion coefficients of a molecule in a cell differ from that in water due to the crowding effect of other molecules. One approach is to model the motion of a molecule as diffusion with constant coefficient, but then have the  effective diffusion be modified by interaction with other particles which are also included in the model (Ridgway, et al., 2008). Another is to 
not include the crowding particles in the model, but to model their effect with a modified diffusion coefficient (Hall and Hoshino, 2010).  Given the inhomogeneity of the cytoplasm, we can expect that this effective diffusion coefficient of the molecule (and its equilibrium probability density)  will vary with location within the cell.
Our framework and numerical method are developed with this latter situation in mind.

In Section~\ref{sec:conclusion} we conclude by explaining the relation between our results and three different 
interpretations of state-dependent diffusion: It\^{o}, Stratonovich, and the Isothermal convention of Lan{\c c}on, Batrouni, Lobry, and Ostrowski (2001).

\section{A Model System for State-Dependent Diffusion}\label{sec:model_sys}

In the Introduction, we described a two-dimensional system consisting of a single particle diffusing inside a rectangular box and reflecting off the boundaries. The diffusion coefficient is twice as large on the right side of the box as the left.  In this section, we demonstrate how to construct a family of deterministic Hamiltonian systems that approximates this behaviour on a coarse scale.

In Subsection~(\ref{subsec:rlg}), we describe the random Lorentz gas (Dettmann, 2000), a deterministic system that when given a random initial condition yields constant-coefficient diffusion at a coarse scale. In Subsection~(\ref{subsec:boxtwo}), we show how to approximate the model system of the introduction by creating two adjacent domains of the random Lorentz gas within a bounding box.  In Subsection~(\ref{subsec:numeric}), we describe numerical experiments with the box system demonstrating that the fraction of time a particle spends on each side of the box cannot be determined solely from  the  values of the diffusion coefficient.
In Subsection~(\ref{subsec:analdisc}), we explain how the result of the numerical experiments in Subsection~(\ref{subsec:numeric}) can be predicted from properties of the dynamics of the microscopic system.


\subsection{Random Lorentz Gas} \label{subsec:rlg}

\begin{figure} 
\begin{center}
\includegraphics[width=6.cm]{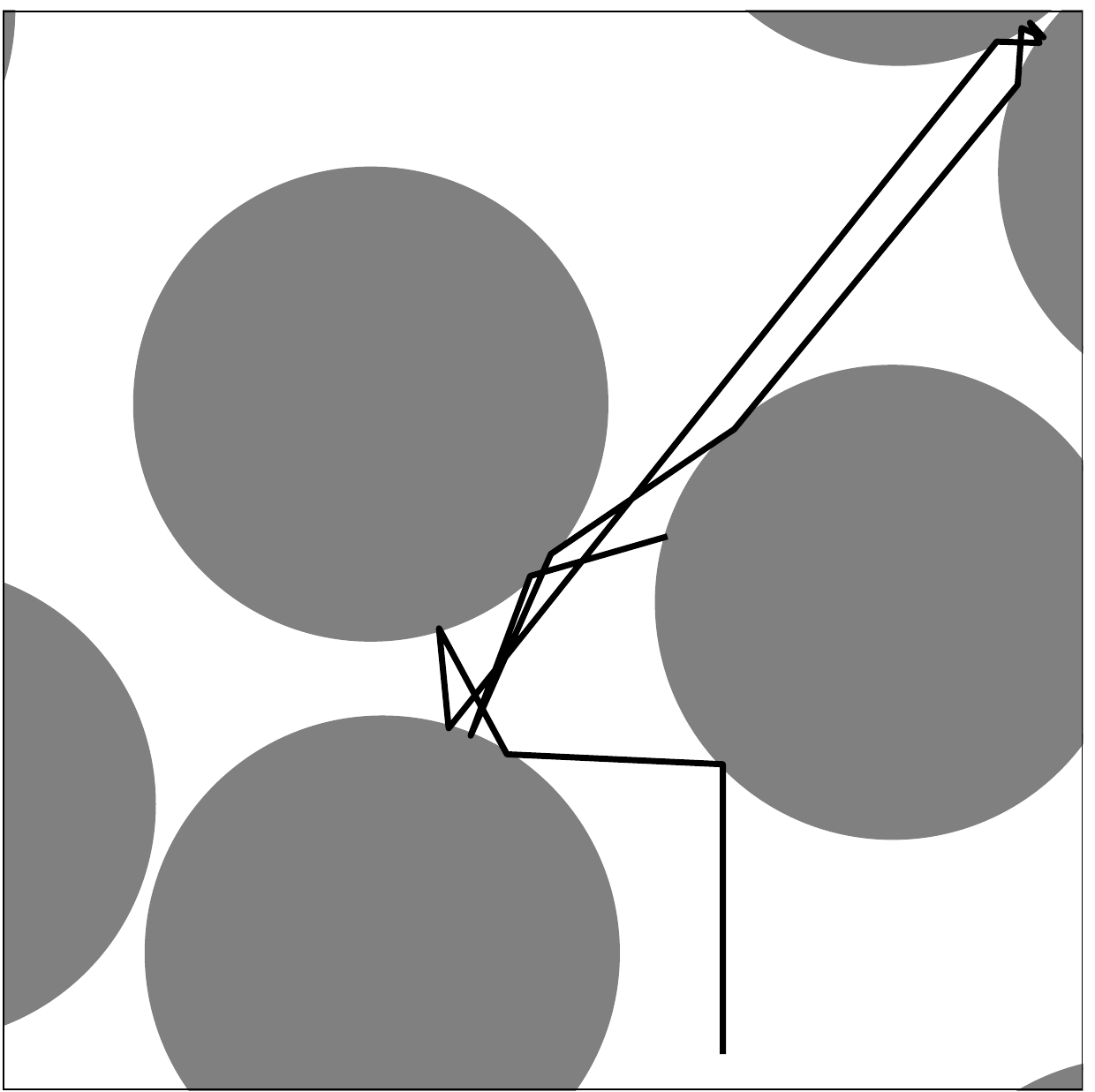}
\hspace{.5cm}
\includegraphics[width=6.cm]{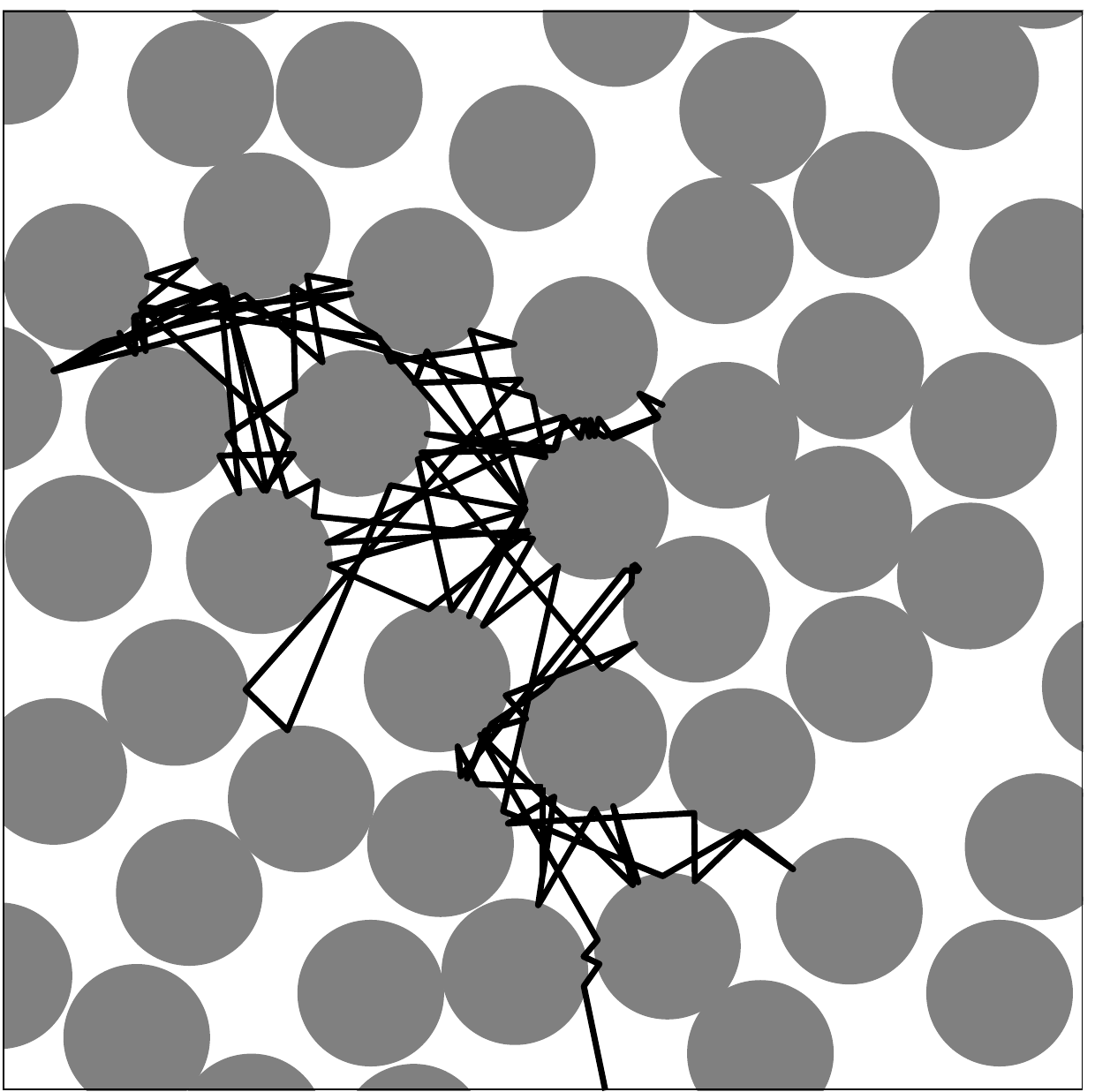}
\end{center}
\caption{ \label{fig:lorentz} Portions of the trajectory of the random Lorentz gas with $\phi=0.5$ at two levels of magnification, showing the discs and the trajectory of the moving particle.}
\end{figure}

Consider infinitely many discs with positions fixed in $\mathbb{R}^2$.
The centres of the discs are distributed randomly with uniform density subject to the constraint that the discs do not overlap.  We consider a single point particle
interacting with the discs in the following manner.  Given an initial velocity and an initial position not on a disc, the  particle moves with constant velocity until it meets a disc.  Then it undergoes an instantaneous elastic reflection with the boundary of the disc, with the angle of incidence equalling the angle of reflection.  This model is  the two-dimensional \emph{random Lorentz gas} (Dettmann, 2000), a mathematical formulation of a model originally due to Lorentz (Lorentz, 1905).  We will always consider the case when the  particle has initial velocity of magnitude 1. Figure~\ref{fig:lorentz} shows trajectories of the random Lorentz gas at two different scales.

The random Lorentz gas has two parameters: the radius of the discs $r$, and the number of discs per unit area $\lambda$.  The radius $r$ can take any positive value.  The density $\lambda$ has  a maximum value $\lambdamax = 1/(r^2 \sqrt{12})$  corresponding to the close-packed hexagonal pattern of discs.  
We define $\phi$ to be the free volume fraction,
the proportion of the area not occupied by discs. 
 We have that $\phi= 1- \pi r^2 \lambda$. The free volume can take any value in $[\phimin, 1)$ where    $\phimin = 1-\pi/\sqrt{12} \approx 0.093$, regardless of the value of $r$.  
 The pair $(r,\phi)$ provides an alternative parametrization of the random Lorentz gas, with the advantage that $\phi$ is dimensionless.

For $\phi = \phimin$ adjacent discs are touching and the  particle remains in a small region of the plane for its entire trajectory. 
For $\phi \in (\phimin,1)$
 the probability that two discs are touching anywhere in the plane is zero (Dettmann 2000, p. 325) and motion of the particle is conjectured to be diffusive (Dettmann and Cohen, 2000). (This contrasts with the situation of overlapping discs for which there is a phase transition to anomalous transport for dense enough placement of discs (H{\"o}fling, Munk, Frey, and Franosch, 2008)).    Specifically, if we start the particle at initial position $x(0)$ not coincident with a disc and give it initial speed 1 and uniformly distributed direction in $[0,2\pi)$, then $x(t)$, the position of the particle at time $t$ is approximately distributed as a Gaussian random vector with mean 0 and variance matrix $2DtI$. Here $I$ is the $2 \times 2$ identity matrix and $D$ denotes the diffusion coefficient.
  This conjecture is supported by analytical calculations (Dettmann, 2000; Ernst and Weyland, 1971)  and numerical simulations (Dettmann and Cohen, 2000; Bruin, 1972).

 A stronger and more formal statement of the conjecture is stated in the language of weak convergence (Billingsley, 1999). Specifically, let $x(t)$ denote the position of the moving particle at time $t$. Fix $x(0)$ to be some point not coincident with a disc and let the initial velocity be chosen as above.
 It is conjectured that 
\[
x(nt)/\sqrt{n} \Rightarrow \sqrt{ 2 D_{r,\phi}} B(t)
\]
as $n \rightarrow \infty$, where $B(t)$ is standard two-dimensional Brownian motion (meaning $\langle B(t)\rangle = 0$, $\langle B(t) B(t)^T \rangle =I$),  $D_{r,\phi}$ is the diffusion coefficient  which depends on $r$ and $\phi$, and $\Rightarrow$ denotes weak convergence in the space of continuous functions (Billingsley, 1999). Such a result holds for the certain periodic Lorentz gasses (Bunimovich and Sina{\u\i}, 1980/81; Klages and Dellago, 2000), but remains open for the random Lorentz gas. (We chose for our study not to use the standard periodic Lorentz gas with discs centred on a hexagonal lattice since for large enough $\phi$ the particle undergoes superdiffusive motion in this case.)

A scaling argument shows that, for fixed $\phi$, $D_{r,\phi}$ is proportional to $r$.  To see this, letting $| \cdot |$ denote the Euclidean norm, observe that the coefficient $D$ can be obtained as $\lim_{t \rightarrow \infty} \langle | x(t) |^2 \rangle / 4t$, assuming $x(0)=0$.  If we rescale space by a factor $R$, we increase both the size of the discs and the distance the particle travels by a factor of $R$, without changing $\phi$. So $\langle | x(t)^2 | \rangle$ increases by a factor of $R^2$.To maintain the speed of the particle as $1$, we also have to rescale time, increasing $t$ by a factor $R$.
 The net effect on the ratio $\langle | x(t) |^2 \rangle / 4 t$ is to increase it by a factor $R$  (Sanders 2005, Sec.\ 3.1.6).
  So
\begin{equation}  \label{eqn:getD}
D_{r,\phi} = r f(\phi)
\end{equation}
for some  function $f$ of $\phi$.

\begin{figure}
\begin{center}
\includegraphics[width=14cm]{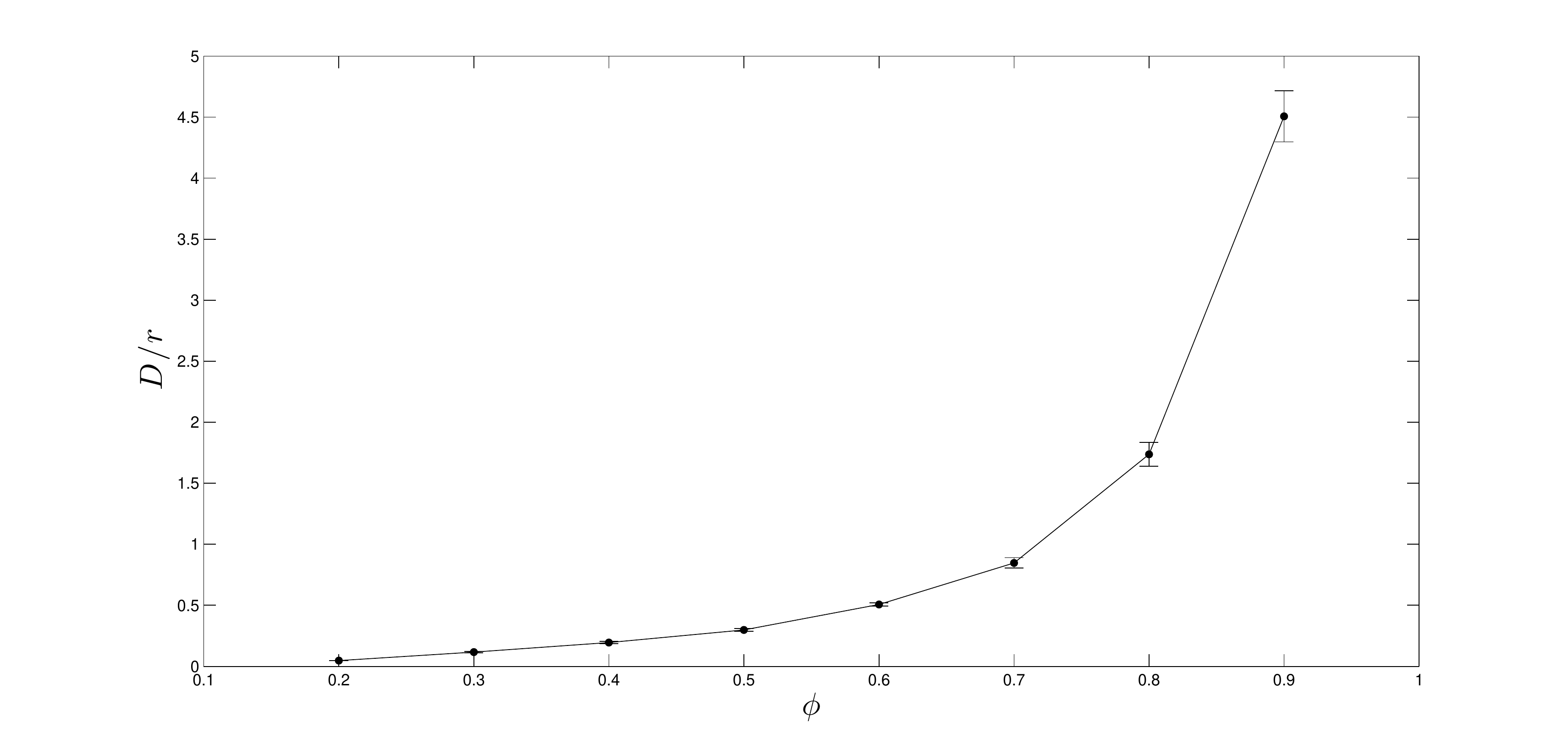}
\end{center}
\caption{ \label{fig:Dvsphi} The ratio of diffusion coefficient versus disc radius $D/r$ as a function of free volume fraction $\phi$ for the random Lorentz gas. }
\end{figure}

Figure~\ref{fig:Dvsphi}  shows the relation between $f(\phi)$ and $\phi$ that was computed using the techniques similar to those described in (Dettmann and Cohen, 2000).
It appears that the function $f(\phi)$ is continuous on its domain, it is monotonically increasing, $f(\phi) \rightarrow 0$ as $\phi \rightarrow \phimin$ and $f(\phi)$ goes to infinity as $\phi \rightarrow 1$. Indeed, calculations from kinetic theory  show that $f(\phi) \sim 3 \pi /[16(1-\phi)]$ in the $\phi \rightarrow 1$ limit 
(vanLeeuwen and Weijland, 1967; Bruin, 1972).

Given any fixed diffusion coefficient $D > 0$ there is a one-parameter family of choices of $r, \phi$ such that $D=D_{r,\phi}$: for any $\phi \in (\phimin, 1)$, just choose $r= D/f(\phi)$.

\subsection{Box with Two Domains} \label{subsec:boxtwo}

In order to investigate the main question of this paper, we take a rectangular box  and divide it into two equal regions.  Each region is filled with randomly placed  discs, with different $r$ and $\phi$ on each side. 
As in the Lorentz gas, the centres of the discs are placed uniformly at random with the condition that they not overlap with each other. Discs can intersect with the walls of the box but not the dividing line between the two sides of the box. Figures \ref{fig:first_params} and \ref{fig:second_params} show two examples of such configurations of discs. The dynamics of the point particle are the same as in the infinite random Lorentz gas, with the added condition that the particle reflects off the walls of the box.

Suppose we are given diffusion coefficients $D_1, D_2>0$.  We can choose $r_1, \phi_1$ and $r_2, \phi_2$ such that $D_1=D_{r_1,\phi_1}$ and $D_2=D_{r_2,\phi_2}$.
 For small enough $r_1, r_2$ the dynamics of the particle in this system will be well-approximated by a particle that diffuses with coefficient $D_1$ on the left side of the box and diffuses with coefficient $D_2$ on the right side of the box.  With appropriate choices of the parameters, we can investigate the question of the proportion of time the particle spends on each side of the box.

\subsection{Numerical Experiments} \label{subsec:numeric}

We start the particle off at some position in the box not on a disc with velocity of magnitude 1 and randomly chosen direction.  The motion of the particle is simulated with an event-driven simulation, computing a trajectory that is accurate up to the errors of floating point arithmetic (Dettmann and Cohen, 2000). 
 Periodically the position of the particle is recorded.  At the end of a long trajectory we compute the number of times the particle  is on the right side of the box divided by the number of times the particle is on the left side of the box.  
 
 We consider two choices of parameters on each side of the box, or set-ups,  each of which leads to the effective diffusion coefficients satisfying $D_2= 2 D_1$.  
 We have contrived the set-ups so that in Set-up 1, the  Statistical Mechanics Prediction is correct, and that in Set-up 2, the  
Time-Change
 Prediction is correct.
 The parameters for each set-up are summarized in Table~\ref{tab:params}
 
 \begin{table}
 \centering
  \begin{tabular}{c|ccc|ccc|c}
Set-up &  $r_1$ & $\phi_1$ & $D_1$ & $r_2$ & $\phi_2$ & $D_2$ & Time on Right / Time on Left  \\
\hline
1 &   0.3   &    0.5  & 0.09 & 0.6    &  0.5  & 0.18 & 0.99 \\
2 &   0.09 &    0.60  & 0.047 &  0.75 & 0.30  & 0.093 & 0.50
\end{tabular}
\caption{ \label{tab:params} The parameters used in each side of the box in the two set-ups, and the ratio between amount of time spent by the particle on the right side of the box and the left side of the box in each set-up.}
\end{table}

\begin{figure} 
\includegraphics[width=14cm]{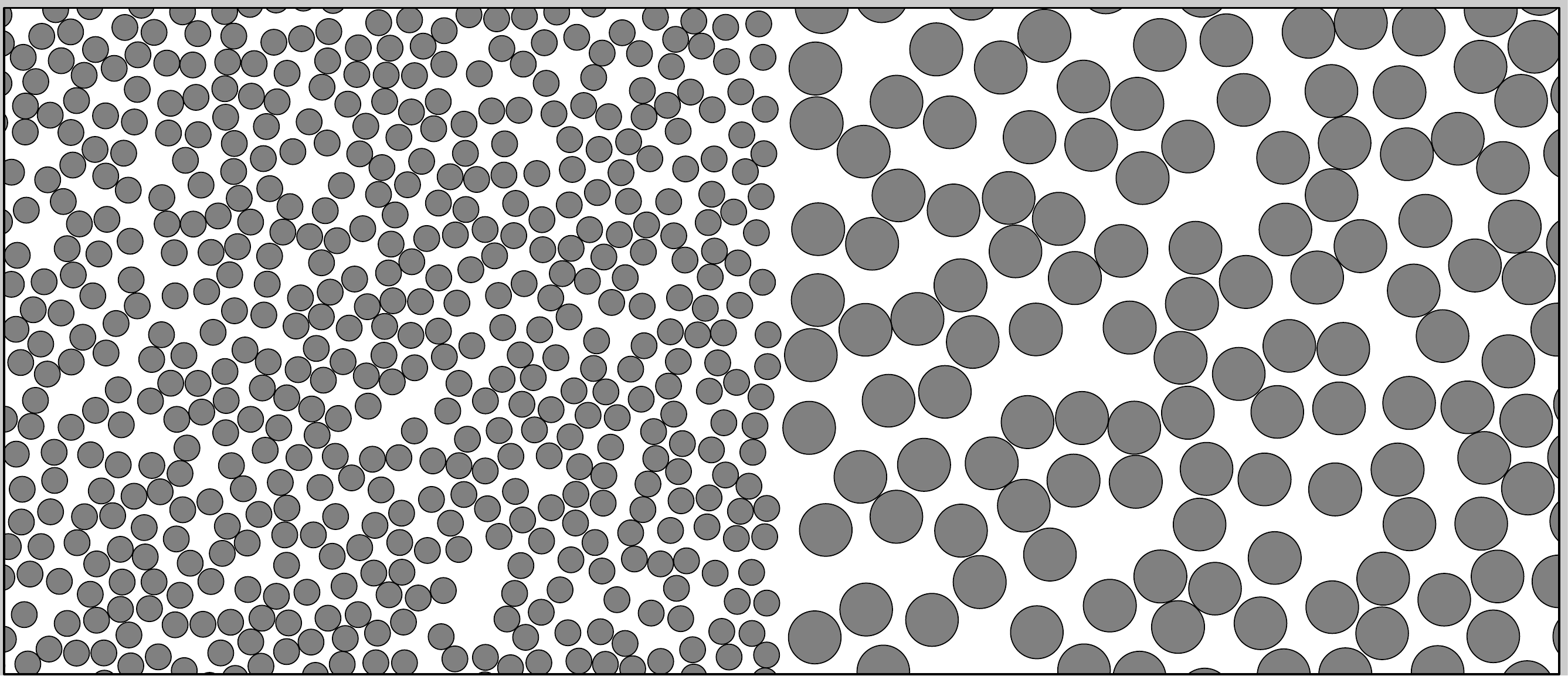}
\caption{\label{fig:first_params} Set-up 1. Free volume fraction is the same on each side: $\phi_1=\phi_2$.  Scatterer radius on the right is twice that on the left: $2r_1=r_2$, leading to $2D_1=D_2$.}
\end{figure}

 In the first set-up $\phi_1= \phi_2 = 0.5 $ and $2 r_1 =  r_2 = 0.6$. We show the position of the discs in the box in Figure~\ref{fig:first_params}.  Over a trajectory of length $5 \times 10^5$ time units the ratio between the occupation times is approximately 1, as we show  in the table above.   This result agrees with the  Statistical Mechanics Prediction.

\begin{figure} 
\includegraphics[width=14cm]{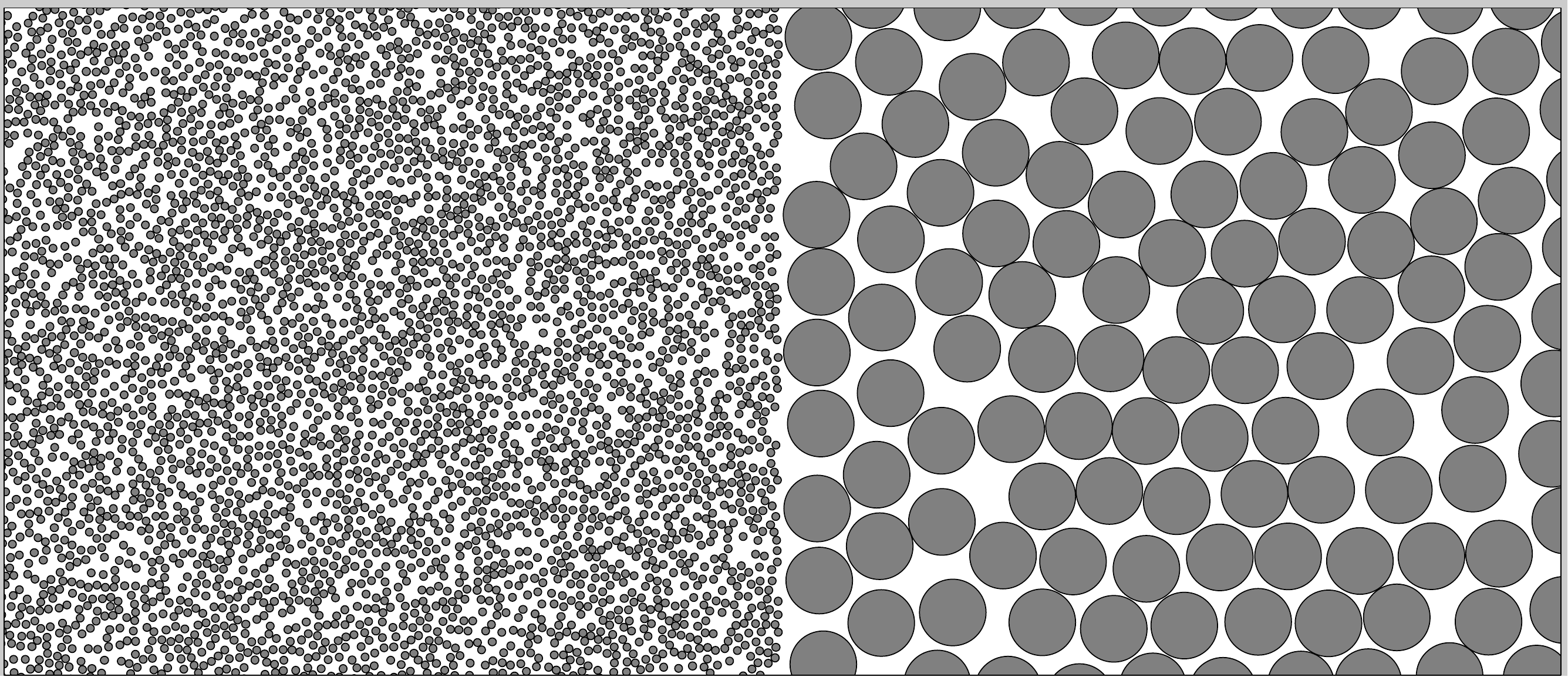}
\caption{ \label{fig:second_params} Set-up 2. Free volume fraction is twice as big on the left: $\phi_1=2\phi_2$. Scatterer radii are chosen so that $2 D_1= D_2$.  }
\end{figure}

 In the second set-up $2 \phi_2 =  \phi_1 =0.60$ and $8.3 r_1 \approx  r_2 = 0.75$. We show the position of the discs in the box in Figure~\ref{fig:second_params}. Over a trajectory of length $7 \times 10^6$ time units the ratio between the occupation times is approximately 1/2, as we show  in the table above. This result agrees with the Time-Change Prediction. 

Thus, by fixing the parameters appropriately, both the Statistical Mechanics Prediction and the
Time-Change Prediction can be seen to be correct for the given mesoscopic behaviour.

\subsection{Analysis and Discussion} \label{subsec:analdisc}


We first explain that, given the properties of random Lorentz gas, the results of the simulation in the previous subsection are predictable.
We first note that the system is ergodic.  The ergodicity of  periodic Lorentz gas is shown by Sina{\u\i} (1970), since it is equivalent to dispersing billiards on the torus. Our system is not dispersing because the walls of the box  are not convex, but Sina{\u\i} (1970, Sec.\ 9) shows how ergodicity still holds in this case by unfolding the box to obtain a periodic domain.

Ergodicity implies that the duration of time that the particle spends in a region of phase space is proportional to the volume of the region.  For our system, this implies that the amount of time spent by the particle on each side is proportional to the free volume fraction $\phi$ on each side.  For fixed $\phi$, the parameters $r$ and $D$ are irrelevant to the proportion of time the particle spends on each side of the box.   In Set-up 1 above, $\phi_1 = \phi_2$ so the ratio of times spent on each side of the box is equal and the Statistical Mechanics Prediction is correct.  In Set-up 2, $\phi_1 = 2 \phi_2$, and so the particle spends twice as much time on the left side of the box, and the 
Time-Change
Prediction is correct.

Despite appearances, the principles of statistical mechanics are not violated in Set-up 2.  There are two ways to reconcile the apparent contradiction.  Firstly, one can say that at the microscopic scale statistical mechanics is not violated because there are not an equal number of states on each side.  The number of states is proportional to the free volume fraction on both sides, and so the system spends more time where the free volume fraction is higher.  The other way of reconciling the disagreement is at the mesoscale. Suppose we divide the box into many 
rectangular cells of equal area, each much smaller than the whole box, but much larger than the size of the discs.
Each cell corresponds to a mesoscopic state which the particle may be in.  Since the system is in the microcanonical ensemble (no exchange of energy with the outside) the probability of the system in equilibrium being in a particular mesoscopic state is determined by the entropy of that state, where the entropy of a mesoscopic state is proportional to the logarithm of the amount of microscopic states it contains. The greater free volume fraction on the left in Set-up 2, implies that mesoscopic states have greater entropy there, and so the system will spend more time on the left side of the box.

On the other hand, in Set-up 1, the Time-Change Prediction proves to be wrong. This means that the motion of the particle in the box is not well-described by the drift-free It\^{o} stochastic differential equation \eqref{eqn:orig_sde}.  Since, by the properties of the uniform random Lorentz gas,  \eqref{eqn:orig_sde} \emph{is} a good model for the dynamics of the particle within each of the two regions of constant disc radius $r$, it must be that the equation is no longer a good model 
at the boundary of the two regions.

We see that there are choices of $\phi_1, r_1$ and $\phi_2, r_2$ such that either  the 
Time-Change Prediction or the Statistical Mechanics Prediction are correct, while still $D_2=2 D_1$.
Indeed, for any $D_1$ and $D_2$, parameters can be chosen to induce arbitrary ratios between the times spent on the left-hand side and the right-hand side.

Although for generic values of the parameters in our Lorentz gas model neither prediction will be valid,  we point out that the Statistical Mechanics Prediction holds for a natural set of parameter settings, whereas  the same is not true of the  Time-Change Prediction.
For the Time-Change 
Prediction to be correct, it is necessary for the ratio between $\phi_1$ and $\phi_2$ to match the ratio between $D_2$ and $D_1$.  We see no natural way that the parameters in our model may be set for this matching to occur. 
On the other hand, for the Statistical Mechanics Prediction to be correct it is necessary that $\phi_1$ and $\phi_2$ be equal. There is at least one case where this condition approximately holds  for a naturally occurring system.
Consider a situation in which both $\phi_1, \phi_2 \approx 1$.  This requires no fine tuning, only that the discs take up a small fraction of the total volume. Choosing $r_1$ and $r_2$ to be unequal leads to different diffusion coefficients on each side, but the particle still spends approximately equal proportions of time in each region.  Likewise, in any physical system where a particle diffuses by interacting with small, sparsely placed scatterers, we expect the Statistical Mechanics Prediction to be correct.

\section{A Proposal for modelling with state-dependent diffusion} \label{sec:prop}

The numerical experiments of the previous section demonstrate that the equilibrium density $\rhoeq(x)$ is not determined solely by the local diffusion rate $D(x)$, even in situations with no external forces acting on the particle.  This raises the practical issue of how to model systems with state-dependent diffusion.   Ideally, mesoscopic diffusive models would be derived from microscopic models via an asymptotic technique such as the van Kampen system-size expansion (van Kampen 2007, Ch.\ XI.3).  But in many circumstances, deriving a realistic microscopic model may be impractical.

Instead, in Subsection~(\ref{subsec:themodel}) we describe a more phenomenological approach.  We assume that the modeller posits an isotropic state-dependent diffusion rate $D(x)$ and an equilibrium density $\rhoeq(x)$.
We then derive a drift coefficient $a(x)$ for an It\^{o} stochastic differential equation with diffusion coefficient $D(x)$ that  gives the desired $\rhoeq(x)$.
In Subsection~(\ref{subsec:connection}) we show how the Lorentz gas system of Section~\ref{sec:model_sys} can be modeled at a mesoscopic level in this way.
In Subsection~(\ref{subsec:themethod}) we then describe an algorithm for simulating such stochastic differential equations (SDEs) that does not refer to $a(x)$ but is expressed in terms only of $D(x)$ and $\rhoeq(x)$.

\subsection{A modelling framework for state-dependent diffusion} \label{subsec:themodel}

We model the diffusion in $k$ spatial dimensions by the It\^{o} SDE
\begin{equation} \label{eqn:framework}
dX(t)  =  a(X) dt + \sqrt{2 D(X)} dB(t), 
\end{equation}
where we will determine $a(x)$ in terms of $\rhoeq(x)$ and $D(x)$. Here $B(t)$ is standard $k$-dimensional Brownian motion.  (In alternate notation, we write this equation as 
$dX(t)/dt = a(x) + \sqrt{2 D(X)} \eta(t)$ where $\eta$ is $k$-dimensional Gaussian white noise.)
 The Fokker-Planck equation for this system is
\[
\frac{\partial}{\partial t} \rho(x,t) = -\nabla \cdot \left[  a(x) \rho(x,t) \right] + 
 \Delta [ D(x) \rho(x,t)] = -\nabla \cdot J(x)
\]
where $J$ is the probability flux.  
Since the problem is underdetermined as stated, we will stipulate that the probability flux vanishes in equilibrium, which is the same as detailed balance holding; see van Kampen (2007, Ch.\ XI.4) for criteria on systems under which this condition holds.
We choose $a(x)$ so that $J(x)$ is zero in equilibrium:
\[
J(x) = a(x) \rhoeq(x)   -  
\nabla [ D(x) \rhoeq(x) ] = 0.
\]
Solving for $a(x)$ gives
\begin{equation} \label{eqn:adef}
a(x) =  \frac{1}{\rhoeq (x) }  \nabla ( D(x) \rhoeq(x))  
=  \nabla D(x) +  D(x) \nabla \ln \rhoeq(x).
\end{equation}
Thus given an equilibrium density $\rhoeq(x)$ and diffusion coefficient $D(x)$ the appropriate It\^{o} SDE is
\begin{equation} \label{eqn:proposed_SDE}
dX(t) = \left( \nabla D(X) +  D(X) \nabla \ln \rhoeq(X) \right) dt + \sqrt{2 D(X)} dB(t).
\end{equation}
The Fokker-Planck equation of \eqref{eqn:proposed_SDE}  is
\begin{equation*}
\frac{\partial}{\partial t} \rho(x,t) 
 =  \nabla \cdot \left[  - \nabla(D(x) \rhoeq(x) ) (\rho(x)/\rhoeq(x))   + \nabla( D(x) \rho(x)) \right],
 \end{equation*}
 or simplifying,
 \begin{equation} \label{eqn:nicepde}
 \frac{\partial}{\partial t}\rho(x,t) =  \nabla \cdot \left[ D(x) \rhoeq(x) \nabla ( \rho(x) /\rhoeq(x) ) \right],
\end{equation}
which is the special isotropic case of Equation XI.4.14 of van Kampen (2007).
Note that if  $\rhoeq(x)$ is constant with respect to $x$, which corresponds to the Statistical Mechanics Prediction being true, then $a(x) = \nabla D(x)$ and \eqref{eqn:proposed_SDE} reduces to 
\[
dX(t) =   \nabla D(x) dt + \sqrt{2 D(x)} dB(t).
\]
On the other hand, to obtain a drift-free It\^{o} SDE ($(a(x)=0$) in this framework requires $D(x) \rhoeq(x)$ to be constant in $x$, which means that $\rhoeq(x)$ is determined completely by $D(x)$.


\subsection{Connection to the Lorentz gas model} \label{subsec:connection}

We explain the connection between the framework of Subsection~(\ref{subsec:themodel}) 
 and the microscopic Lorentz gas model of Section~\ref{sec:model_sys}. Given an instance of our random Lorentz gas model with two domains, what are the corresponding functions $D(x)$, $\rhoeq(x)$ in  \eqref{eqn:nicepde}, the mesoscopic equation for $\rho$?
 
 Suppose on the left side of the box the Lorentz gas model has parameters $r_1$ and $\phi_1$ and on the right it has parameters $r_2$ and $\phi_2$. The diffusion coefficients on each side, $D_1$ and $D_2$ are determined by the relation \eqref{eqn:getD}.  To determine the equilibrium probability density of the particle on each side, let $2A$ be the total area of the box, so that each side has area $A$. We know that $\rhoeqi$, the probability density on side $i$, is proportional to $\phi_i$, and thus $\rhoeqone/\rhoeqtwo = \phi_1/\phi_2$. We also know that since the total probability must be $1$, $\rhoeqone A + \rhoeqtwo A = 1$. Solving for $\rhoeqi$ gives
\[
\rhoeqi = \frac{\phi_i}{A(\phi_1 + \phi_2)},
\]
for $i=1,2$.
We define  $D(x)$ and $\rhoeq(x)$ for $x=(x_1,x_2)$ in the box  by
\[
D(x) = \left\{ 
\begin{array}{ll} 
D_1, & \mbox{for }x_1 < 0, \\ 
D_2, & \mbox{for }x_1 >0,
\end{array}\right.
\hspace{1cm}
\rhoeq(x) = \left\{ 
\begin{array}{ll} 
\rhoeqone, & \mbox{for }x_1 < 0, \\ 
\rhoeqtwo, & \mbox{for }x_1 >0.
\end{array} \right.
\]
These functions determine $a(x)$, the drift coefficient in \eqref{eqn:framework}, via 
\eqref{eqn:adef}.  The drift $a(x)$ is zero everywhere except along the line $x_1=0$ where it is not defined. 

We determine the appropriate boundary conditions for $\rho$ along the boundary line $x_1=0$. 
In order for the right hand side of  \eqref{eqn:nicepde} to be well defined we require $\rho(x)/\rhoeq(x)$ to be continuous. So for any $x$ along the boundary line we must have
\[
\frac{\rho(x^-)}{\rhoeq(x^-)} = \frac{\rho(x^+)}{\rhoeq(x^+)}
\]
where $x^-$ denotes taking the limit from the left, and $x^+$ denotes taking the limit from the right.
For our particular choice of $\rhoeq$ this gives
\[
\frac{\rho(x^-)}{\rhoeqone} = \frac{\rho(x^+)}{\rhoeqtwo}.
\]
The second boundary condition comes from assuming continuous flux across the boundary line:
\[
D(x^-) \rhoeq(x^-) \frac{\partial}{\partial x_1} \frac{\rho(x^-)}{\rhoeq(x^-)} =
D(x^+) \rhoeq(x^+) \frac{\partial}{\partial x_1} \frac{\rho(x^+)}{\rhoeq(x^+)}.
\] 
For our particular choices of $D$ and $\rhoeq$, since $\rhoeq$ is constant away from the line $x_1=0$, this gives the boundary conditions
\[
D_1 \frac{\partial}{\partial x_1} \rho(x^-) = D_2 \frac{\partial}{\partial x_1} \rho(x^+).
\]

A possible direction for further investigation is to consider Lorentz gas models where disc radius $r$ and free volume fraction $\phi$ vary smoothly with $x$, and to determine what $D(x)$ and $\rhoeq(x)$, and hence $a(x)$ are in this case.

\subsection{Numerical Simulation} \label{subsec:themethod}

If both $\rhoeq$ and $D$ are smooth  then the Euler-Maruyama scheme or Milstein's method  is effective for simulating \eqref{eqn:proposed_SDE}  (see Higham 2001, for example).  If either $D(x)$ or $\rhoeq(x)$ are discontinuous with respect to $x$, then the drift term will have a singularity that is not resolved by  standard time-integration schemes, and numerical simulations may not correctly approximate the equilibrium density.   We propose using a variant of the  Metropolis-adjusted Langevin algorithm (MALA) as introduced by Roberts and Tweedie (1996) and analyzed by Bou-Rabee and Vanden-Eijnden (2010).   MALA is obtained by taking a convergent numerical method for the SDE (in this case, Euler-Maruyama) and introducing Metropolis step-rejections in order to have a discrete-time process with the correct equilibrium density. Our approach is  different in that we start with a convergent method for only the diffusive part of the It\^{o} SDE, and then we introduce  step rejections to induce the correct drift and equilibrium density.

Let $h$ be the step length of our numerical discretisation and let $X_n$ be the numerical approximation to $X(nh)$, where $X$ is  the  solution to the It\^{o} SDE \eqref{eqn:proposed_SDE} for $n \geq 0$.
Given a numerical value $X_n$, we let the trial value for the next step be defined by 
\begin{equation}  \label{eqn:eulermaru}
X^*_{n+1} = X_n + \sqrt{2 D(X_n)} [ B((n+1)h) - B(nh) ]
\end{equation}
where $B$ is standard $k$-dimensional Brownian motion.
Then  $X_{n+1}$  is given by
\begin{equation} \label{eqn:metropolis}
X_{n+1} = \left\{ \begin{array}{ll}
X^*_{n+1}  &  \mbox{    if } \xi_k <  \alpha_h(X_n,X^*_{n+1}),    \\
X_n & \mbox{    otherwise,}
\end{array}  \right.
\end{equation}
where 
\[
\alpha_h(x,y) = \min \left(1, \frac{q_h(y,x) \rhoeq( y)  }{q_h(x,y) \rhoeq(x) }\right).
\]
and $\xi_k, k \geq 1$ is an independent, identically distributed sequence of random variables, uniform on $[0,1]$ and independent of $B$.
Here 
\[
q_h(x,y) = \frac{1}{\sqrt{4 \pi h D(x)}} e^{- (x-y)^2 /4 h D(x)}
\]
is the transition probability density for $X^*_{n+1}$ being at $y$ given that $X_n$ is at $x$.
The definition of $X^*_{n+1}$ in \eqref{eqn:eulermaru} is just the Euler-Maruyama method for the It\^{o} stochastic differential equation $dX = \sqrt{2 D(X)} dB$.  
The Metropolis rejection procedure \eqref{eqn:metropolis} accepts the Euler-Maruyama step with probability $\alpha_h(X_n,X^*_{n+1})$ and rejects it otherwise.

The Metropolis rejection procedure guarantees that the process $X_n, n\geq 1$ has $\rhoeq(x)$ as its equilibrium density.  On the other hand, in any region of the state space where $D(x)$ and $\rhoeq(x)$ are constant with respect to $x$, $\alpha_h$ is $1$ and so the method reduces to the Euler-Maruyama method for the constant coefficient diffusion without drift. Future work will investigate rigorously the convergence of the above method to the solutions of \eqref{eqn:proposed_SDE}.

Here we numerically demonstrate the convergence of the method for the simple case where 
\[
D(x) = \left\{ \begin{array}{ll}
1 & \mbox{   if } x \leq 0, \\
2 & \mbox{   if } x > 0, 
\end{array} \right.
\]
and 
\[
\rhoeq(x) = \left\{ \begin{array}{ll}
1 & \mbox{   if } x \in [-1,1], \\
0 & \mbox{  otherwise.}
\end{array} \right.
\]
Equation $\eqref{eqn:proposed_SDE}$ with this choice of $D(x)$ and $\rhoeq(x)$ models a particle  diffusing on the interval $[-1,1]$ with reflecting boundary conditions at $\pm 1$ and a piecewise constant diffusion coefficient.  The reflecting boundary conditions are conveniently implemented by our choice of $\rhoeq(x)$. 

\begin{figure}
\includegraphics[width=6.8cm]{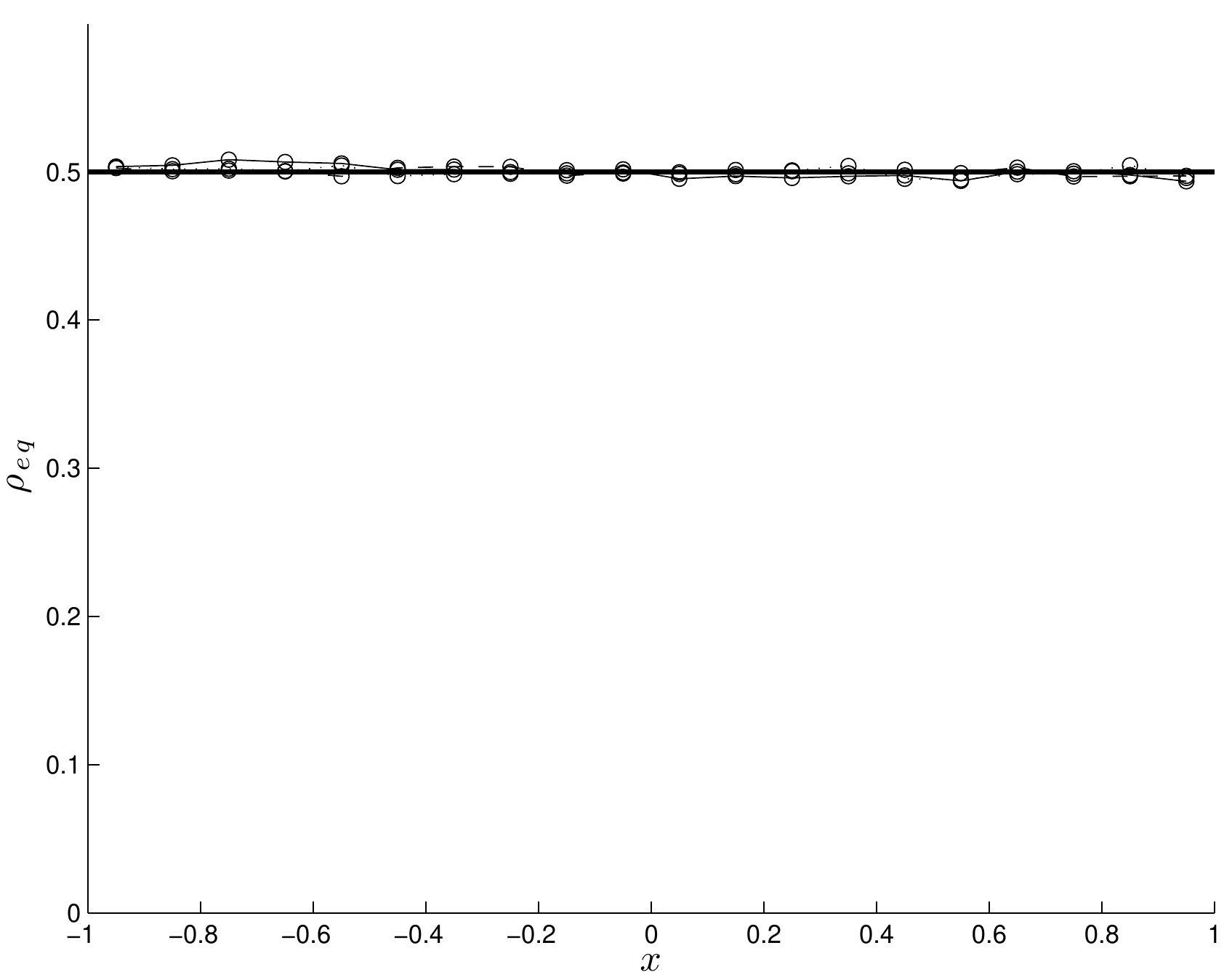}
\includegraphics[width=6.8cm]{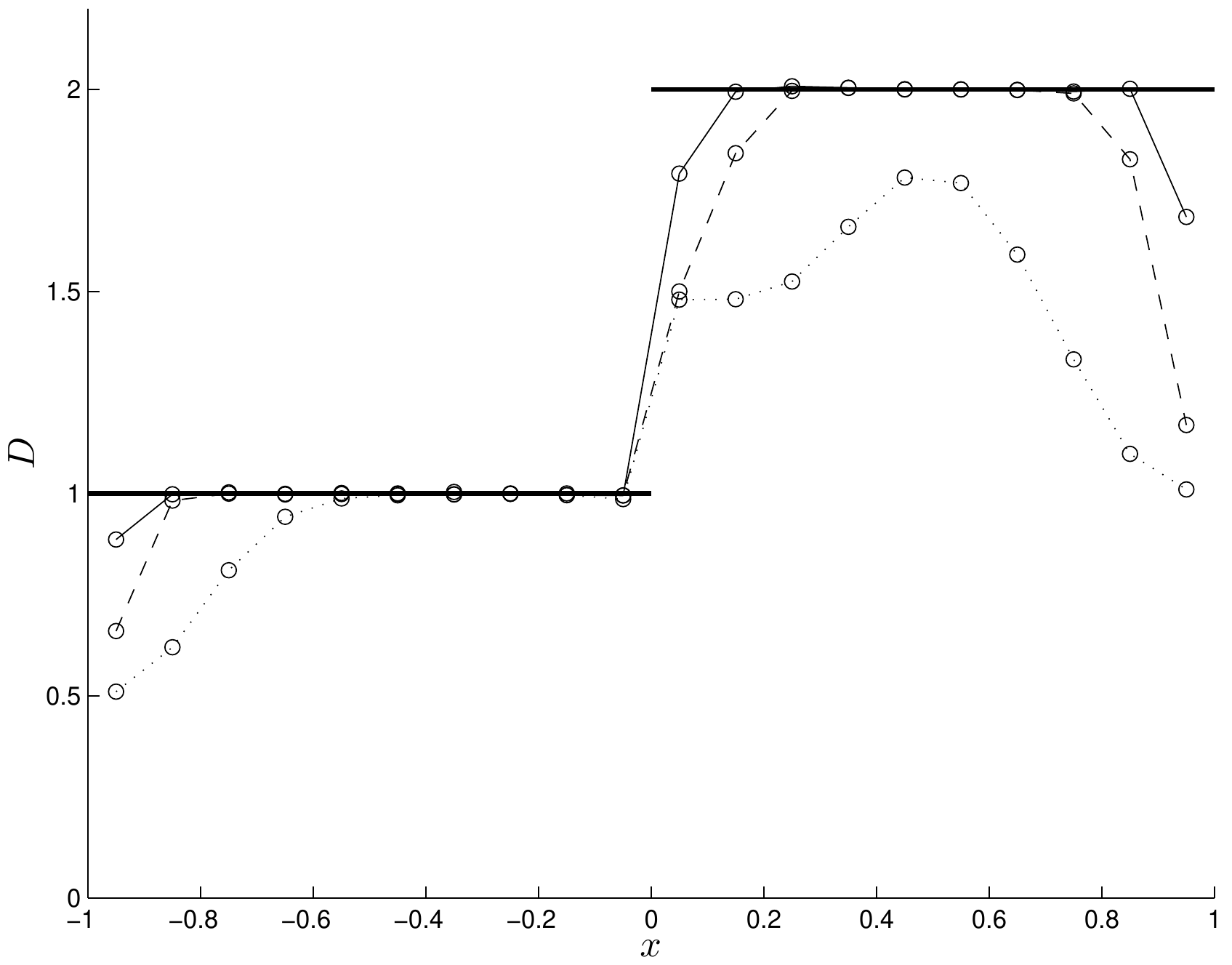}
\caption{\label{fig:newmethod}  Equilibrium density $\rhoeq$ and local diffusion coefficient $D$ for method \eqref{eqn:eulermaru} applied to a simple one-dimensional SDE.  Shown are results for $h=0.01$ (dotted), $h=0.001$ (dashed), and $h=0.0001$ (solid), along with the exact values for the SDE.}
\end{figure}

Figure~\ref{fig:newmethod} shows the results of simulating \eqref{eqn:proposed_SDE} with these choices of $D(x)$ and $\rhoeq(x)$ using the method we have described. To show results, we divide the interval $[-1,1]$ into 20 equal subintervals, and plot the density for the amount of time the particle spends in each subinterval. We also plot the effective diffusion coefficient for each bin, which we define to the the average observed value of $(X_{n+1}-X_n)^2/2h$ over the trajectory, for all $n$ such that $X_n$ lies in the given bin. Simulations were conducted with  $h=0.01, 0.001, 0.0001$ and for trajectories long enough so that standard statistical errors in the plots are smaller than the symbols used.
We see that for all values of $h$ the equilibrium density $\rhoeq(x)$ is correctly reproduced. The effective diffusion coefficient converges to $D$ as $h$ goes to zero.

\section{Discussion: It\^{o}, Stratonovich, or Isothermal} \label{sec:conclusion}

We conclude by discussing our results in the  context of the apparent ambiguity between It\^{o}, Stratonovich, and Isothermal interpretations of stochastic integrals (Lau and Lubensky, 2007; Volpe, et al., 2010; Kupferman, et al., 2004).
For the purposes of discussion, we consider a particle moving in one spatial dimension whose position at time $t$ is $X(t)$. We assume that $X(t)$ is a Markov stochastic process with continuous sample paths.  We model the motion of the particle with the stochastic differential equation (SDE)
\begin{equation} \label{eqn:generalSDE}
dX(t) = a(X(t)) dt + b(X(t)) dB
\end{equation}
where $B(t)$ is standard Brownian motion.
We call $a(x)$  the drift and  $b(x)$ the diffusion of the SDE. 

As is well known (van Kampen, 2007; Gardiner, 2004), unless we specify a particular interpretation, the SDE \eqref{eqn:generalSDE}  does not unambiguously define the stochastic process $X(t)$.
To see this, we integrate \eqref{eqn:generalSDE} over $[0,T]$ to get
\[
X(t)-X(0) = \int^T_0 a(X(t)) dt + \int_0^T b(X(t)) dB(t).
\]
The first term on the right has a unique interpretation as a Riemann integral, but the second term cannot be simply viewed as a Riemann-Stieljes integral, since $B$ is not of bounded variation.
If we compute the second term as the limit of Riemann sums, the answer depends on where in each subinterval the argument $b(X(t))$ is evaluated.
For example, choosing $h=T/N$, $t_n=hn$ and letting $B_n = B(t_n)$ and $X_n = X(t_n)$,
suppose we take the integral with respect to $B$ to be
\[
\int_0^T b(X(t)) dB(t) = \lim_{h \rightarrow 0} \sum_{n=0}^{N-1}  b(X_n^*) (B_{n+1} -B_n),
\]
where 
\[
X_n^* = (1-\alpha) X_n + \alpha X_{n+1}.
\]
Famously, unless $b(x)$ is a constant, the limit depends on the choice of $\alpha$ (Volpe, et al., 2010). If we choose $\alpha=0$, we obtain the It\^{o} interpretation of the integral, which yields a stochastic process $X(t)$ with Fokker-Planck equation
\begin{equation*} \label{eqn:ito-fokker-planck}
\frac{\partial}{\partial t} \rho(x,t) = -\frac{\partial}{\partial x} \left[ a(x) \rho(x,t) \right] + \frac{1}{2} \frac{\partial^2}{\partial x^2} \left[ b(x)^2 \rho(x,t) \right].
\end{equation*}
 If we choose  $\alpha=1/2$, we obtain the Stratonovich interpretation of the integral, which yields a process $X(t)$ with Fokker-Planck equation
 \begin{eqnarray*} \label{eqn:strat-fokker-planck}
\frac{\partial}{\partial t} \rho(x,t) 
& = &  -\frac{\partial}{\partial x} \left[ a(x) \rho(x,t) + \frac{1}{2} b(x) b'(x) \rho(x,t) \right] +    \frac{1}{2} \frac{\partial^2}{\partial x^2} \left[ b(x)^2 \rho(x,t) \right]. 
\end{eqnarray*}
 Note that  the Fokker-Planck equation shows that the Stratonovich interpretation of the SDE with drift $a(x)$ and diffusion $b(x)$ yields the same stochastic process as the It\^{o} interpretation of the SDE with drift $a(x) + b(x) b'(x) /2$ and diffusion $b(x)$ (Gardiner 2004, p.\ 99).
 Finally, if we choose $\alpha=1$ we obtain the  anti-It\^{o} or Isothermal interpretation (Lau and Lubensky, 2007; Volpe, et al., 2010), which yields a process with Fokker-Planck equation 
  \begin{eqnarray*} \label{eqn:iso-fokker-planck}
\frac{\partial}{\partial t} \rho(x,t) 
& = &  -\frac{\partial}{\partial x} \left[ a(x) \rho(x,t) +  b(x) b'(x) \rho(x,t) \right] +    \frac{1}{2} \frac{\partial^2}{\partial x^2} \left[ b(x)^2 \rho(x,t) \right]. 
\end{eqnarray*}
In this case, the Isothermal interpretation of the SDE with drift $a(x)$ and diffusion $b(x)$ gives the same stochastic process as the It\^{o} interpretation of the SDE with drift $a(x)+b(x)b'(x)$ and diffusion $b(x)$ (Lau and Lubensky, 2007).

As we can see in the various Fokker-Planck equations above, if we fix $a(x)$ and $b(x)$, varying the parameter $\alpha$ gives different stochastic processes for  the motion of the particle. This fact may make it seem like there should be a physically correct choice of the parameter $\alpha$.  We argue that this is false. For any fixed $\alpha$ the range of stochastic processes that can be captured by an appropriate choice of $a(x)$ and $b(x)$ is the same.  For example, suppose we fix a choice of $a(x)$ and $b(x)$ and choose to interpret \eqref{eqn:generalSDE} with a given $\alpha \in [0,1]$.  The process defined is identical to what we would obtain with the It\^{o} interpretation ($\alpha=0$) of the SDE with  drift  $a(x)+ \alpha b'(x) b(x)$ and diffusion  $b(x)$. 

Though the families of stochastic processes described using each convention are the same, it may still be the case that some choice of $\alpha$ is more natural or convenient for some purposes than others. Frequently, the rationale is based on the idea that if the drift is zero, then $X(t)$ should have certain properties.
For example, if one wants $\langle X(t) -X(0)\rangle = 0$ for all $t$ when $a(x) \equiv 0$, regardless of $b(x)$, then the It\^{o} convention with $\alpha=0$ guarantees this.
 If one wants that when $a(x)\equiv 0$ the equilibrium density is constant, then the Isothermal convention with $\alpha=1$ guarantees this. 
We summarize the properties of the various interpretations of the SDE \eqref{eqn:generalSDE} when $a(x) \equiv 0$  in Table~\ref{tab:interp}.

\begin{table}
\centering
\begin{tabular}{c|c|c|c}
$\alpha$ & \parbox{3.5cm}{Interpretation of Stochastic Integral} & \parbox{3.5cm}{ Fokker-Plank Equation for $a \equiv 0$}  & \parbox{3.5cm}{Special Properties when $a \equiv 0$} \\ \hline
        0 &   It\^{o}   &  
        $\frac{\partial }{\partial t} \rho = \frac{1}{2} \frac{\partial^2}{\partial x^2} \left[ b^2 \rho \right]$ &  
        $\langle X(t) - X(0) \rangle =0$ \mbox{ for all $t$.} 
        \\
        1/2 & Stratonovich & 
          $\frac{\partial }{\partial t} \rho = \frac{1}{2} \frac{\partial}{\partial x} \left\{ b \frac{\partial}{\partial x} \left[ b \rho \right]\right\}$ &
   \\
        1 &   Isothermal   &  
        $\frac{\partial }{\partial t} \rho = \frac{1}{2} \frac{\partial}{\partial x} \left\{ b^2 \frac{\partial}{\partial x} \rho\right\}$ &
         $\rhoeq = $ const.
        \\
\end{tabular}
\caption{ \label{tab:interp} Summary of some properties of the It\^{o}, Stratonovich, and Isothermal interpretations of \eqref{eqn:generalSDE}.}
\end{table}

The question we posed in the Introduction may be rephrased as follows: in the absence of external forces, and given a diffusion $b(x)= \sqrt{2D(x)}$, what is the correct choice of parameter $\alpha$ and drift $a(x)$ to model the motion of the particle?  A natural way to approach the problem is to interpret the absence of external forces as meaning that $a(x)\equiv 0$.  Then the problem boils down to the choice of $\alpha$: the Statistical-Mechanics Prediction follows from taking $\alpha=1$, and the Time-Change Prediction follows from taking $\alpha=0$.  The results in Section~\ref{sec:model_sys} showed that neither answer is justified universally. 

In Section~\ref{sec:prop} we recommended a different approach.  We fix $\alpha$ and then choose $a(x)$ to generate the desired equilibrium distribution.  As we have explained here, the choice of $\alpha$ is not crucial once we allow a non-zero $a(x)$.  Accordingly, we have chosen $\alpha=0$, corresponding to It\^{o} calculus. This is the main choice in the mathematics literature, and numerical methods such as the Euler-Maruyama method take a particularly simple form with it.  Once we have made this choice of $\alpha$, we are free to choose $a(x)$ appropriately. In Section~\ref{sec:prop} we chose $a(x)$ to ensure a given equilibrium density.
 
Beyond the particular needs of the present work, we believe the framework we describe in Section~\ref{sec:prop} provides a natural and flexible way to model diffusive systems. 
  In situations where a researcher is confident for physical reasons that the equilibrium probability is constant, then our framework takes a simple form. 
  The numerical method we present functions even when the diffusion and equilibrium density are discontinuous. Future work will study the convergence properties of the method, as well exploring applications of our general framework.

\ack{The authors would like to thank John Bechhoefer for suggesting the problem,  and Florian Theil, Michael Plischke, and Martin Zuckermann for helpful discussions. We also thank Nilima Nigam, John Bechhoefer, and Nancy Forde for extensive comments on an earlier draft of this manuscript.}


{\bf References}

P.~Billingsley.
 {\em Convergence of probability measures}.
 John Wiley \& Sons Inc., 2nd edition, 1999.

N.~Bou-Rabee and E.~Vanden-Eijnden.
 Pathwise accuracy and ergodicity of metropolized integrators for
  {SDE}s.
 {\em Communications on Pure and Applied Mathematics}, 63(5):655--696,
  2010.

C.~Bruin.
 Logarithmic terms in the diffusion coefficient for the lorentz gas.
 {\em Phys. Rev. Lett.}, 29:1670--1674, Dec 1972.

L.~A. Bunimovich and Y.~G. Sina{\u\i}.
 Statistical properties of {L}orentz gas with periodic configuration
  of scatterers.
 {\em Comm. Math. Phys.}, 78(4):479--497, 1980/81.

C.~P. Dettmann.
 The {L}orentz gas: a paradigm for nonequilibrium stationary states.
 In {\em Hard ball systems and the {L}orentz gas}, volume 101 of {\em
  Encyclopaedia Math. Sci.}, pages 315--365. Springer, Berlin, 2000.

C.~P. Dettmann and E.~G.~D. Cohen.
 Microscopic chaos and diffusion.
 {\em J. Statist. Phys.}, 101(3-4):775--817, 2000.

M.~Ernst and A.~Weyland.
 Long time behaviour of the velocity auto-correlation function in a
  {L}orentz gas.
 {\em Phys. Lett. A}, 34(1):39--40, January 1971.

C.~W. Gardiner.
 {\em Handbook of Stochastic Methods for Physics, Chemistry and the
  Natural Sciences}.
 Springer, 3rd edition, 2004.

D.~Hall and M.~Hoshino.
{\em Effects of macromolecular crowding on intracellular diffusion from a single particle perspective.}
{\em Biophys. Rev.}, 2(1):39--53, 2010.

D.~J. Higham.
 An algorithmic introduction to numerical simulation of stochastic
  differential equations.
 {\em SIAM Review}, 43(3):525--546, 2001.

F.~H{\"o}fling, T.~Munk, E.~Frey, and T.~Franosch.
 Critical dynamics of ballistic and brownian particles in a
  heterogeneous environment.
 {\em The Journal of Chemical Physics}, 128(16):164517, 2008.

R.~Klages and C.~Dellago.
 Density-dependent diffusion in the periodic lorentz gas.
 {\em J. Stat. Phys.}, 101(1--2):145--159, Oct 2000.

N.~Korabel and E.~Barkai.
 Boundary conditions of normal and anomalous diffusion from thermal
  equilibrium.
 {\em Phys. Rev. E}, 83:051113, 2011.

R.~Kupferman, G~.A.~Pavliotis, and A.~M.~Stuart.
It\^o versus Stratonovich white-noise limits for systems with inertia and colored multiplicative noise.
{Phys. Rev. E}, 70:036120, 2004.

P.~Lan{\c c}on, G.~Batrouni, L.~Lobry, and N.~Ostrowsky.
 Drift without flux: Brownian walker with a space-dependent diffusion
  coefficient.
 {\em Europhysics Letters}, 54:28--34, 2001.

A.~W.~C. Lau and T.~C. Lubensky.
 State-dependent diffusion: Thermodynamic consistency and its path
  integral formulation.
 {\em Phys. Rev. E}, 76(1):011123, Jul 2007.

H.~A. Lorentz.
 The motion of electrons in metallic bodies.
 {\em Proc. Roy. Acad. Amst.}, 7:438--453, 1905.

H.~G.~Othmer and A.~Stevens.
Aggregation, blowup, and collapse: The ABC's of taxis in reinforced random walks.
{\em SIAM J. Appl. Math.} 57(4):1044--1081, 1997.

O.~Ovaskainen and S.~J. Cornell.
 Biased movement at a boundary and conditional occupancy times for
  diffusion processes.
 {\em J. Appl. Probab.}, 40(3):557--580, 2003.

F.~Reif.
 {\em Fundamentals of Statistical and Thermal Physics}.
 McGraw-Hill, New York, 1965.

D.~Ridgway, G.~Broderick, A.~Lopez-Campistrous, M.~Ru'aini, P. Winter, M.~Hamilton, P.~Boulanger, A.~Kovalenko, and M.~J.~Ellison.
Coarse-grained molecular simulation of diffusion and reaction kinetics in a crowded virtual cytoplasm.
{\em Biophys. J.},  94(10):3748--3759, 2008. 

G.~O. Roberts and R.~L. Tweedie.
 Exponential convergence of {L}angevin distributions and their
  discrete approximations.
 {\em Bernoulli}, 2(4):341--363, 1996.

D.~Sanders.
 {\em Deterministic Diffusion in Periodic Billiard Models}.
 PhD thesis, University of Warwick, 2005.

J.~G. Sina{\u\i}.
 Dynamical systems with elastic reflections. {E}rgodic properties of
  dispersing billiards.
 {\em Uspehi Mat. Nauk}, 25(2 (152)):141--192, 1970.

T.~E. Turner, S.~Schnell, and K.~Burrage.
Stochastic approaches for modelling in vivo reactions.
{\em Comp Bio and Chem}, 28(3):165--178, 2004.

N.~G. van Kampen.
 {\em Stochastic Processes in Physics and Chemistry}.
 North-Holland Personal Library. North-Holland, 3rd edition, 2007.

J.~vanLeeuwen and A.~Weijland.
 Non-analytic density behaviour of the diffusion coefficient of a
  lorentz gas i. divergencies in the expansion in powers in the density.
 {\em Physica}, 36(3):457--490, 1967.

G.~Volpe, L.~Helden, T.~Brettschneider, J.~Wehr, and C.~Bechinger.
 Influence of noise on force measurements.
 {\em Phys. Rev. Lett.}, 104(17):170602, 2010.


\end{document}